\documentclass[12pt,reqno]{article}

\usepackage[usenames]{color}
\usepackage{amssymb}
\usepackage{graphicx}
\usepackage{amscd}

\usepackage[colorlinks=true,
linkcolor=webgreen,
filecolor=webbrown,
citecolor=webgreen]{hyperref}

\definecolor{webgreen}{rgb}{0,.5,0}
\definecolor{webbrown}{rgb}{.6,0,0}

\usepackage{color}
\usepackage{fullpage}
\usepackage{float}

\usepackage{graphics,amsmath,amssymb}
\usepackage{amsfonts}
\usepackage{latexsym}
\usepackage{epsf}

\newcommand{\seqnum}[1]{\href{http://www.research.att.com/cgi-bin/access.cgi/as/~njas/sequences/eisA.cgi?Anum=#1}{\underline{#1}}}

\begin{document}

\begin{center}
\vskip 1cm{\LARGE\bf 
An efficient algorithm for the computation of Bernoulli numbers \\
\vskip .1in
\vskip 1cm
\large}
Greg Fee\\
Centre for Constructive and Experimental Mathematics \\
Simon Fraser University\\
Vancouver, 
Canada\\
\href{mailto:gfee@cecm.sfu.ca}{\tt gfee@cecm.sfu.ca} \\
\ \\
Simon Plouffe \\
Montr\'eal, 
Canada\\
\href{mailto:simon.plouffe@gmail.com}{\tt mailto:simon.plouffe@gmail.com} \\
\end{center}

\vskip .2 in
\begin{abstract}
This article gives a direct formula for the computation of $B \left( n \right)$ using the
asymptotic formula $$B \left( n \right) \approx 2\,{\frac {n!}{{\pi }^{n}{2}^{n}}}$$
where $n$ is even and $n \gg 1$.
This is simply based on the fact that $\zeta  \left( n \right)$ is very near 1 when $n$ is large and
since 
$B \left( n \right) = 2\,{\frac {\zeta  \left( n \right) n!}{{
\pi }^{n}{2}^{n}}}$
exactly. The formula chosen for the Zeta function is the one with prime numbers from the well-known Euler product for $\zeta  \left( n \right)$. This algorithm is far better than the recurrence formula for the Bernoulli numbers even
if each $B(n)$ is computed  individually. The author could compute $B \left( 750,000 \right)$ in a few hours. The current record of computation is now (as of Feb. 2007) $B \left( 5,000,000 \right)$ a number
of (the numerator) of 27332507 decimal digits is also based on that idea.\\
\end{abstract}

\newtheorem{theorem}{Theorem}[section]
\newtheorem{proposition}{Proposition}[section]
\newtheorem{corollary}{Corollary}[section]
\newtheorem{lemma}{Lemma}[section]

%
%
%

\newcommand{\ben}{\begin{enumerate}}
\newcommand{\een}{\end{enumerate}}
\newcommand{\ble}{\begin{lem}}
\newcommand{\ele}{\end{lem}}
\newcommand{\bth}{\begin{thm}}
\renewcommand{\eth}{\end{thm}}
\newcommand{\bpr}{\begin{prop}}
\newcommand{\epr}{\end{prop}}
\newcommand{\bco}{\begin{cor}}
\newcommand{\eco}{\end{cor}}
\newcommand{\bcon}{\begin{conj}}
\newcommand{\econ}{\end{conj}}
\newcommand{\bde}{\begin{defn}}
\newcommand{\ede}{\end{defn}}
\newcommand{\bex}{\begin{exa}}
\newcommand{\eex}{\end{exa}}
\newcommand{\barr}{\begin{array}}
\newcommand{\earr}{\end{array}}
\newcommand{\btab}{\begin{tabular}}
\newcommand{\etab}{\end{tabular}}
\newcommand{\beq}{\begin{equation}}
\newcommand{\eeq}{\end{equation}}
\newcommand{\bea}{\begin{eqnarray*}}
\newcommand{\eea}{\end{eqnarray*}}
\newcommand{\bce}{\begin{center}}
\newcommand{\ece}{\end{center}}
\newcommand{\bpi}{\begin{picture}}
\newcommand{\epi}{\end{picture}}
\newcommand{\bfi}{\begin{figure} \begin{center}}
\newcommand{\efi}{\end{center} \end{figure}}
\newcommand{\capt}{\caption}
\newcommand{\bsl}{\begin{slide}{}}
\newcommand{\esl}{\end{slide}}

\newcommand{\bib}{thebibliography}
\newcommand{\pf}{{\bf Proof}\hspace{7pt}}
\newcommand{\qed}{\rule{1ex}{1ex}}
\newcommand{\Qed}{\rule{1ex}{1ex} \medskip}
\newcommand{\qqed}{\qquad\rule{1ex}{1ex}}
\newcommand{\Qqed}{\qquad\rule{1ex}{1ex}\medskip}
\newcommand{\mc}[3]{\multicolumn{#1}{#2}{#3}}
\newcommand{\Mc}[1]{\multicolumn{#1}{c}{}}
\newcommand{\ul}{\underline}
\newcommand{\ol}{\overline}
\newcommand{\hor}{\mbox{--}}
\newcommand{\hs}[1]{\hspace{#1}}
\newcommand{\hso}[1]{\hspace{-1pt}}
\newcommand{\vs}[1]{\vspace{#1}}
\newcommand{\qmq}[1]{\quad\mbox{#1}\quad}
\newcommand{\ssk}{\smallskip}
\newcommand{\msk}{\mediumskip}
\newcommand{\bsk}{\bigskip}
\newcommand{\rp}[2]{\rule{#1pt}{#2pt}}
\newcommand{\st}[1]{\rule{#1pt}{0pt}}
\newcommand{\mbl}[1]{\makebox(0,0)[l]{#1}}
\newcommand{\mbr}[1]{\makebox(0,0)[r]{#1}}
\newcommand{\mbt}[1]{\makebox(0,0)[t]{#1}}
\newcommand{\mbb}[1]{\makebox(0,0)[b]{#1}}
\newcommand{\mbc}[1]{\makebox(0,0){#1}}

\newcommand{\emp}{\emptyset}
\newcommand{\indu}{\uparrow}
\newcommand{\rest}{\downarrow}
\newcommand{\sm}{\hspace{-2pt}\setminus\hspace{-2pt}}
\newcommand{\sbs}{\subset}
\newcommand{\sbe}{\subseteq}
\newcommand{\sps}{\supset}
\newcommand{\spe}{\supseteq}
\newcommand{\setm}{\setminus}
\newcommand{\asy}{\thicksim}
\newcommand{\dle}{<\hspace{-6pt}\cdot}
\newcommand{\dge}{\cdot\hspace{-6pt}>}
\newcommand{\iso}{\cong}
\newcommand{\Cong}{\equiv}
\newcommand{\conm}[1]{\stackrel{#1}{\equiv}}

\newcommand{\widevec}{\overrightarrow}
\newcommand{\Ch}{\hat{C}}
\newcommand{\eh}{\hat{e}}
\newcommand{\Gh}{\hat{G}}
\newcommand{\ih}{\hat{\imath}}
\newcommand{\jh}{\hat{\jmath}}
\newcommand{\gh}{\hat{g}}
\newcommand{\mh}{\hat{0}}
\newcommand{\Mh}{\hat{1}}
\newcommand{\Dh}{\hat{D}}
\newcommand{\zh}{\hat{0}}
\newcommand{\oh}{\hat{1}}
\newcommand{\Ph}{\hat{P}}
\newcommand{\ph}{\hat{p}}
\newcommand{\Pih}{\hat{\Pi}}
\newcommand{\Qh}{\hat{Q}}
\newcommand{\Sh}{\hat{S}}
\newcommand{\lh}{\hat{l}}
\newcommand{\xh}{\hat{x}}
\newcommand{\ptn}{\vdash}
\newcommand{\lt}{\lhd}
\newcommand{\gt}{\rhd}
\newcommand{\lte}{\unlhd}
\newcommand{\gte}{\unrhd}
\newcommand{\jn}{\vee}
\newcommand{\Jn}{\bigvee}
\newcommand{\mt}{\wedge}
\newcommand{\Mt}{\bigwedge}
\newcommand{\bdy}{\partial}
\newcommand{\sd}{\bigtriangleup}
\newcommand{\case}[4]{\left\{\barr{ll}#1&\mbox{#2}\\#3&\mbox{#4}\earr\right.}
\newcommand{\fl}[1]{\lfloor #1 \rfloor}
\newcommand{\ce}[1]{\lceil #1 \rceil}
\newcommand{\flf}[2]{\left\lfloor\frac{#1}{#2}\right\rfloor}
\newcommand{\cef}[2]{\left\lceil\frac{#1}{#2}\right\rceil}
\newcommand{\gau}[2]{\left[ \barr{c} #1 \\ #2 \earr \right]}
\newcommand{\setc}[2]{\{left{#1}\ :\ #2\}}
\newcommand{\setl}[2]{\{left{#1}\ |\ #2\}}
\newcommand{\qst}{$q$-Stirling numbers}
\newcommand{\qbi}{$q$-binomial coefficients}
\newcommand{\del}{\nabla}
\newcommand{\pde}[2]{\frac{\partial#1}{\partial#2}}
\newcommand{\spd}[2]{\frac{\partial^2\hspace{-2pt}#1}{\partial#2^2}}
\newcommand{\ode}[2]{\ds\frac{d#1}{d#2}}
\def\<{\langle}
\def\>{\rangle}
\newcommand{\fall}[2]{\langle{#1}\rangle_{#2}}
\newcommand{\nor}[1]{\parallel{#1}\parallel}
\newcommand{\ipr}[2]{\langle{#1},{#2}\rangle}
\newcommand{\spn}[1]{\langle{#1}\rangle}
\newcommand{\ree}[1]{(\ref{#1})}
\newcommand{\rpl}{$\longleftarrow$}

\newcommand{\ra}{\rightarrow}
\newcommand{\Ra}{\Rightarrow}
\newcommand{\lta}{\leftarrow}
\newcommand{\rla}{\leftrightarrow}
\newcommand{\lra}{\longrightarrow}
\newcommand{\lla}{\longleftarrow}
\newcommand{\Lra}{\Longrightarrow}
\newcommand{\llra}{\longleftrightarrow}
\newcommand{\Llra}{\Longleftrightarrow}
\newcommand{\ve}[1]{\stackrel{\ra}{#1}}

\newcommand{\al}{\alpha}
\newcommand{\be}{\beta}
\newcommand{\ga}{\gamma}
\newcommand{\de}{\delta}
\newcommand{\ep}{\epsilon}
\newcommand{\io}{\iota}
\newcommand{\ka}{\kappa}
\newcommand{\la}{\lambda}
\newcommand{\mut}{\tilde{\mu}}
\newcommand{\om}{\omega}
\newcommand{\Phit}{\tilde{\Phi}}
\newcommand{\Psit}{\tilde{\Psi}}
\newcommand{\rhot}{\tilde{\rho}}
\newcommand{\si}{\sigma}
\renewcommand{\th}{\theta}
\newcommand{\ze}{\zeta}
\newcommand{\Al}{\Alpha}
\newcommand{\Ga}{\Gamma}
\newcommand{\De}{\Delta}
\newcommand{\Ep}{\epsilon}
\newcommand{\La}{\Lambda}
\newcommand{\Om}{\Omega}
\newcommand{\Si}{\Sigma}
\newcommand{\Th}{\Theta}
\newcommand{\Ze}{\Zeta}

\newcommand{\0}{{\bf 0}}
\newcommand{\1}{{\bf 1}}
\newcommand{\2}{{\bf 2}}
\newcommand{\3}{{\bf 3}}
\newcommand{\4}{{\bf 4}}
\newcommand{\5}{{\bf 5}}
\newcommand{\6}{{\bf 6}}
\newcommand{\7}{{\bf 7}}
\newcommand{\8}{{\bf 8}}
\newcommand{\9}{{\bf 9}}
\newcommand{\ba}{{\bf a}}
\newcommand{\bb}{{\bf b}}
\newcommand{\bc}{{\bf c}}
\newcommand{\bd}{{\bf d}}
\newcommand{\bfe}{{\bf e}}
\newcommand{\bff}{{\bf f}}
\newcommand{\bg}{{\bf g}}
\newcommand{\bh}{{\bf h}}
\newcommand{\bi}{{\bf i}}
\newcommand{\bj}{{\bf j}}
\newcommand{\bk}{{\bf k}}
\newcommand{\bm}{{\bf m}}
\newcommand{\bn}{{\bf n}}
\newcommand{\bp}{{\bf p}}
\newcommand{\bq}{{\bf q}}
\newcommand{\br}{{\bf r}}
\newcommand{\bs}{{\bf s}}
\newcommand{\bt}{{\bf t}}
\newcommand{\bu}{{\bf u}}
\newcommand{\bv}{{\bf v}}
\newcommand{\bw}{{\bf w}}
\newcommand{\bx}{{\bf x}}
\newcommand{\by}{{\bf y}}
\newcommand{\bz}{{\bf z}}
\newcommand{\bA}{{\bf A}}
\newcommand{\bB}{{\bf B}}
\newcommand{\bC}{{\bf C}}
\newcommand{\bE}{{\bf E}}
\newcommand{\bF}{{\bf F}}
\newcommand{\bG}{{\bf G}}
\newcommand{\bH}{{\bf H}}
\newcommand{\bN}{{\bf N}}
\newcommand{\bP}{{\bf P}}
\newcommand{\bQ}{{\bf Q}}
\newcommand{\bR}{{\bf R}}
\newcommand{\bS}{{\bf S}}
\newcommand{\bT}{{\bf T}}
\newcommand{\bX}{{\bf X}}
\newcommand{\bXh}{\hat{{\bf X}}}
\newcommand{\bZ}{{\bf Z}}
\newcommand{\bep}{{\bf \ep}}
\newcommand{\bcdot}{{\bf \cdot}}
\newcommand{\bbC}{{\mathbb C}}
\newcommand{\bbF}{{\mathbb F}}
\newcommand{\bbN}{{\mathbb N}}
\newcommand{\bbP}{{\mathbb P}}
\newcommand{\bbQ}{{\mathbb Q}}
\newcommand{\bbR}{{\mathbb R}}
\newcommand{\bbS}{{\mathbb S}}
\newcommand{\bbZ}{{\mathbb Z}}
\newcommand{\cA}{{\cal A}}
\newcommand{\cAB}{{\cal AB}}
\newcommand{\cAD}{{\cal AD}}
\newcommand{\cB}{{\cal B}}
\newcommand{\cBC}{{\cal BC}}
\newcommand{\cC}{{\cal C}}
\newcommand{\cD}{{\cal D}}
\newcommand{\cDB}{{\cal DB}}
\newcommand{\cE}{{\cal E}}
\newcommand{\cF}{{\cal F}}
\newcommand{\cG}{{\cal G}}
\newcommand{\cH}{{\cal H}}
\newcommand{\cI}{{\cal I}}
\newcommand{\cIF}{{\cal IF}}
\newcommand{\cJ}{{\cal J}}
\newcommand{\cK}{{\cal K}}
\newcommand{\cL}{{\cal L}}
\newcommand{\cm}{{\cal m}}
\newcommand{\cM}{{\cal M}}
\newcommand{\cN}{{\cal N}}
\newcommand{\cNBC}{{\cal NBC}}
\newcommand{\cO}{{\cal O}}
\newcommand{\cP}{{\cal P}}
\newcommand{\cPD}{{\cal PD}}
\newcommand{\cR}{{\cal R}}
\newcommand{\cRF}{{\cal RF}}
\newcommand{\cS}{{\cal S}}
\newcommand{\cT}{{\cal T}}
\newcommand{\cV}{{\cal V}}
\newcommand{\cW}{{\cal W}}
\newcommand{\cX}{{\cal X}}
\newcommand{\cY}{{\cal Y}}
\newcommand{\cZ}{{\cal Z}}

\newcommand{\df}{\dot{f}}
\newcommand{\dF}{\dot{F}}
\newcommand{\fgl}{{\mathfrak gl}}
\newcommand{\fS}{{\mathfrak S}}
\newcommand{\at}{\tilde{a}}
\newcommand{\ct}{\tilde{c}}
\newcommand{\dt}{\tilde{d}}
\newcommand{\pt}{\tilde{p}}
\newcommand{\Bt}{\tilde{B}}
\newcommand{\Gt}{\tilde{G}}
\newcommand{\Ht}{\tilde{H}}
\newcommand{\Kt}{\tilde{K}}
\newcommand{\Nt}{\tilde{N}}
\newcommand{\St}{\tilde{S}}
\newcommand{\Xt}{\tilde{X}}
\newcommand{\alt}{\tilde{\alpha}}
\newcommand{\bet}{\tilde{\beta}}
\newcommand{\rht}{\tilde{\rho}}
\newcommand{\Pit}{\tilde{\Pi}}
\newcommand{\tal}{\tilde{\alpha}}
\newcommand{\tbe}{\tilde{\beta}}
\newcommand{\tPi}{\tilde{\Pi}}
\newcommand{\ab}{\ol{a}}
\newcommand{\Ab}{\ol{A}}
\newcommand{\Bb}{\ol{B}}
\newcommand{\cb}{\ol{c}}
\newcommand{\Cb}{{\ol{C}}}
\newcommand{\db}{\ol{d}}
\newcommand{\Db}{\ol{D}}
\newcommand{\eb}{\ol{e}}
\newcommand{\Eb}{\ol{E}}
\newcommand{\fb}{\ol{f}}
\newcommand{\Gb}{\ol{G}}
\newcommand{\Hb}{\ol{H}}
\newcommand{\ib}{\ol{\imath}}
\newcommand{\Ib}{\ol{I}}
\newcommand{\jb}{\ol{\jmath}}
\newcommand{\Jb}{\ol{J}}
\newcommand{\kb}{\ol{k}}
\newcommand{\Kb}{\ol{K}}
\newcommand{\nb}{\ol{n}}
\newcommand{\pb}{\ol{p}}
\newcommand{\Pb}{\ol{P}}
\newcommand{\Phib}{\ol{\Phi}}
\newcommand{\Qb}{\ol{Q}}
\newcommand{\rb}{\ol{r}}
\newcommand{\Rb}{\ol{R}}
\newcommand{\Sb}{\ol{S}}
\newcommand{\tb}{\ol{t}}
\newcommand{\Tb}{\ol{T}}
\newcommand{\Ub}{\ol{U}}
\newcommand{\Wb}{\ol{W}}
\newcommand{\xb}{\ol{x}}
\newcommand{\Xb}{\ol{X}}
\newcommand{\yb}{\ol{y}}
\newcommand{\Yb}{\ol{Y}}
\newcommand{\zb}{\ol{z}}
\newcommand{\Zb}{\ol{Z}}
\newcommand{\pib}{\ol{\pi}}
\newcommand{\sib}{\ol{\si}}
\newcommand{\degb}{\ol{\deg}}
\newcommand{\vj}{\vec{\jmath}}
\newcommand{\vv}{\vec{v}}
\newcommand{\ttab}{\{t\}}
\newcommand{\stab}{\{s\}}

\newcommand{\Aut}{\mathop{\rm Aut}\nolimits}
\newcommand{\abl}{\mathop{\rm al}\nolimits}
\newcommand{\capa}{\mathop{\rm cap}\nolimits}
\newcommand{\codim}{\mathop{\rm codim}\nolimits}
\newcommand{\ch}{\mathop{\rm ch}\nolimits}
\newcommand{\col}{\mathop{\rm col}\nolimits}
\newcommand{\ctt}{\mathop{\rm ct}\nolimits}
\newcommand{\Der}{\mathop{\rm Der}\nolimits}
\newcommand{\des}{\mathop{\rm des}\nolimits}
\newcommand{\Des}{\mathop{\rm Des}\nolimits}
\newcommand{\diag}{\mathop{\rm diag}\nolimits}
\newcommand{\Diag}{\mathop{\rm Diag}\nolimits}
\newcommand{\diam}{\mathop{\rm diam}\nolimits}
\newcommand{\End}{\mathop{\rm End}\nolimits}
\newcommand{\ev}{\mathop{\rm ev}\nolimits}
\newcommand{\eval}{\mathop{\rm eval}\nolimits}
\newcommand{\Eval}{\mathop{\rm Eval}\nolimits}
\newcommand{\FPF}{\mathop{\rm FPF}\nolimits}
\newcommand{\gen}{\mathop{\rm gen}\nolimits}
\newcommand{\Harm}{\mathop{\rm Harm}\nolimits}
\newcommand{\Hom}{\mathop{\rm Hom}\nolimits}
\newcommand{\id}{\mathop{\rm id}\nolimits}
\newcommand{\im}{\mathop{\rm im}\nolimits}
\newcommand{\imaj}{\mathop{\rm imaj}\nolimits}
\newcommand{\inc}{\mathop{\rm inc}\nolimits}
\newcommand{\Ind}{\mathop{\rm Ind}\nolimits}
\newcommand{\inter}{\mathop{\rm int}\nolimits}
\newcommand{\Inter}{\mathop{\rm Int}\nolimits}
\newcommand{\inv}{\mathop{\rm inv}\nolimits}
\newcommand{\Inv}{\mathop{\rm Inv}\nolimits}
\newcommand{\lcm}{\mathop{\rm lcm}\nolimits}
\newcommand{\LHS}{\mathop{\rm LHS}\nolimits}
\newcommand{\Mat}{\mathop{\rm Mat}\nolimits}
\newcommand{\maj}{\mathop{\rm maj}\nolimits}
\newcommand{\Mod}{\mathop{\rm mod}\nolimits}
\newcommand{\NBB}{\mathop{\rm NBB}\nolimits}
\newcommand{\nin}{\mathop{\rm nin}\nolimits}
\newcommand{\Nin}{\mathop{\rm Nin}\nolimits}
\newcommand{\nul}{\mathop{\rm nul}\nolimits}
\newcommand{\od}{\mathop{\rm od}\nolimits}
\newcommand{\per}{\mathop{\rm per}\nolimits}
\newcommand{\Pete}{\mathop{\rm Pete}\nolimits}
\newcommand{\pfa}{\mathop{\rm pf}\nolimits}
\newcommand{\Pm}{\mathop{\rm pm}\nolimits}
\newcommand{\PM}{\mathop{\rm PM}\nolimits}
\newcommand{\Pol}{\mathop{\rm Pol}\nolimits}
\newcommand{\rad}{\mathop{\rm rad}\nolimits}
\newcommand{\RHS}{\mathop{\rm RHS}\nolimits}
\newcommand{\rk}{\mathop{\rm rk}\nolimits}
\newcommand{\rank}{\mathop{\rm rank}\nolimits}
\newcommand{\reg}{\mathop{\rm reg}\nolimits}
\newcommand{\row}{\mathop{\rm row}\nolimits}
\newcommand{\sdiam}{\mathop{\rm sdiam}\nolimits}
\newcommand{\sgn}{\mathop{\rm sgn}\nolimits}
\newcommand{\sign}{\mathop{\rm sign}\nolimits}
\newcommand{\sh}{\mathop{\rm sh}\nolimits}
\newcommand{\slack}{\mathop{\rm slack}\nolimits}
\newcommand{\spa}{\mathop{\rm span}\nolimits}
\newcommand{\sta}{\mathop{\rm st}\nolimits}
\newcommand{\summ}{\mathop{\rm sum}\nolimits}
\newcommand{\sus}{\mathop{\rm susp}\nolimits}
\newcommand{\Sym}{\mathop{\rm Sym}\nolimits}
\newcommand{\SYT}{\mathop{\rm SYT}\nolimits}
\newcommand{\tang}{\mathop{\rm tang}\nolimits}
\newcommand{\tr}{\mathop{\rm tr}\nolimits}
\newcommand{\tri}{\mathop{\rm tri}\nolimits}
\newcommand{\wed}{\mathop{\rm wedge}\nolimits}
\newcommand{\wt}{\mathop{\rm wt}\nolimits}

\newcommand{\bul}{\bullet}
\newcommand{\dil}{\displaystyle}
\newcommand{\dsum}{\dil\sum}
\newcommand{\dint}{\dil\int}
\newcommand{\scl}{\scriptstyle}
\newcommand{\ssl}{\scriptscriptstyle}
\newcommand{\foz}{\footnotesize}
\newcommand{\scz}{\scriptsize}
\newcommand{\cho}{\choose}

\newcommand{\Schu}{Sch\"utzenberger}
\newcommand{\aim}{Adv.\ in Math.\/}
\newcommand{\bams}{Bull.\ Amer.\ Math.\ Soc.\/}
\newcommand{\cjm}{Canad.\ J. Math.\/}
\newcommand{\dm}{Discrete Math.\/}
\newcommand{\dmj}{Duke Math.\ J.\/}
\newcommand{\ejc}{European J. Combin.\/}
\newcommand{\jaa}{J. Algebra\/}
\newcommand{\jac}{J. Algebraic Combin.\/}
\newcommand{\jas}{J. Algorithms\/}
\newcommand{\jams}{J. Amer.\ Math.\ Soc.\/}
\newcommand{\jct}{J. Combin.\ Theory\/}
\newcommand{\jcta}{J. Combin.\ Theory Ser. A\/}
\newcommand{\jctb}{J. Combin.\ Theory Ser. B\/}
\newcommand{\jgt}{J. Graph Theory\/}
\newcommand{\jram}{J. Reine Angew.\ Math.\/}
\newcommand{\pjm}{Pacific J. Math.\/}
\newcommand{\pams}{Proc.\ Amer.\ Math.\ Soc.\/}
\newcommand{\plms}{Proc.\ London Math.\ Soc.\/}
\newcommand{\tams}{Trans.\ Amer.\ Math.\ Soc.\/}
\newcommand{\pja}{Proc.\ Japan Acad.\ Ser.\ A  Math\/}
\newcommand{\sv}{Springer-Verlag Lecture Notes in Math.\/}
\newcommand{\crgs}{Combinatoire et Repr\'{e}sentation du 
    Groupe Sym\'{e}trique, Strasbourg 1976, D. Foata ed.}
\newcommand{\caup}{Cambridge University Press}
\newcommand{\oup}{Oxford University Press} 
\newcommand{\pr}{preprint}
\newcommand{\ip}{in preparation}
\newcommand{\ds}{\displaystyle}

\setlength{\topmargin}{.1in}
\setlength{\textheight}{8in}
\setlength{\textwidth}{7in}
\setlength{\evensidemargin}{-.2in}
\setlength{\oddsidemargin}{-.2in}

\newtheorem{thm}{Theorem}[section]
\newtheorem{prop}[thm]{Proposition}
\newtheorem{cor}[thm]{Corollary}
\newtheorem{lem}[thm]{Lemma}
\newtheorem{conj}[thm]{Conjecture}
\newtheorem{exa}[thm]{Example}
\newtheorem{question}[thm]{Question}

\section{The need for a single computation}

This algorithm came once in 1996 when the authors
wanted to compute large Bernoulli numbers using a well-known computer algebra system system like Maple or
Mathematica. These programs used Faulhaber's recurrence \cite{A&S, hardy} formula which is nice but unsuitable for large computations. We quickly came to the conclusion 
that $B(10000)$ was out of reach even with a powerful
computer. This is where we realized that for $n$ large the
actual formula is simply $B \left( n \right) =2\,{\frac {n!}{{\pi }^{n}{2}^{n}}}$ where n is even and not counting the sign, for $n$=1000 the approximation is good to more than 300 decimal digits where $B \left( 1000 \right)$ is of the order of 1770 digits. To carry out the exact computation of $B \left( 1000 \right)$ one has only to compute first the principal term in the asymptotic formula and secondly just a few terms in the Euler product (up to $p = 59$). 
The second idea was that the fractional part of the Bernoulli numbers can also be computed very fast with the help of the Von Staudt-Clausen formula. So finally, the need is only to compute $B_n$ with 
enough precision so that the remainder is $< 1$ and apply the Von Staudt-Clausen formula for
the fractional part to finally add the 2 results. Note : Mathematica now uses a much more efficient algorithm partly due to these results presented here.

\section{The Von Staudt-Clausen formula}

The formula is, for $k \geq 1$, $$\left( -1 \right) ^{k}B_{2k} \equiv{\it \sum} \left( {\frac {1}{p}} \right)\mod 1.$$
The sum being extended over primes p such that $(p-1)|2k$ \cite {hardy}. In other words, for $B(10)$ the sum is

$$B(10) = 1 - 1/2 - 1/3 -1/11 = 5/66.$$ In terms of computation, when $n$ is
of the order of 1000000 it goes very fast to compute the fractional part of $B_n$.
The only thing that remains to be done then is the principal part or integer part
of $B_n$.

\section{The Euler product}

The Euler product of the zeta function is 

$$\sum_{n=1}^\infty \frac{1}{n^s}= \prod _{p\in\mathbb{P}} \frac{1}{1-p^{-s}}.$$

Where $s > 1$ and $p$ is prime. This is the error term in $B_n$. For any given $n$ there are ${\frac {n}{\ln  \left( n \right) }}$ primes compared to $n$. Translated into the program it means less operations to carry, the program stops when $p^k$ is of the order of $B \left( n \right)$.  

\vspace{2ex} 
\section{The final program}

The Maple program uses a high precision value of $2\pi$ and a routine for the 
Von Staudt-Clausen formula. That program held the record of the computation of 
Bernoulli Numbers from 1996 to 2002, after that others made more efficient programs 
using C++ and high precision packages like Kellner and Pavlyk (see table 1) 
and could reach $B(5,000,000)$.\\
The program was used in 2003 to verify Agoh's conjecture up to n=49999 by the authors.
Agoh's conjecture is $$pB_{p-1} \equiv -1 \mod p$$ is
true iff $p$ is prime. The congruence is not obvious since $pB_{p-1}$ is a fraction. The standard method
reduces first the numerator mod p, then re-evaluates the fraction, then reduces the numerator
mod p. The final fraction is always smaller than 1 and the result of $a/b \mod p$ is solved
by finding $k$ such that $a \equiv bk \mod p$.
There are 3 parts in the main program which may take time. 
First the computation of $({2\pi})^n$ and $n!$. 
Secondly, the evaluation of the Von Staudt-Clausen formula and 
thirdly the computation of the Euler product.\
On a medium sized computer (Pentium 2.4 Ghz with Maple 10 and 1 gigabyte of memory). The run
time for $B(20000)$ is about 9 seconds and the number is 61382 digits long including 1 second
to read the value of $\pi$ to high-precision from the disk.
Here are the timings
for that run :\\

- Product with primes up to 1181 at 61382 digits of precision : 7 seconds.

- Exponentiation of ${2\pi}$ and $n!$ : less than 1 second.

- Computation of 20000! : negligible.

- Computation of Von Staudt-Clausen expression : negligible.\\
When $n$ increases the time taken to evaluate the product with primes is what takes
the most.\
A value of $\pi$ to several thousands digits is necessary. Maple can supply many
thousands but a file containing 1 million is easily found on the internet and is much faster. 
In this program $\pi$ is renamed pi with no capitals.\
The Bernoulli numbers up to $n=100$ are within the program mainly for speed when $n$ is small.\\

\begin{verbatim}

BERN:=proc(n::integer)
local d, z, oz, i, p, pn, pn1, f, s, p1, t1, t2;
global Digits;
    lprint(`start at time` = time());
    if n = 1 then -1/2
    elif n = 0 then 1
    elif n < 0 then ERROR(`argument must be >= 0`)
    elif irem(n, 2) = 1 then 0
    elif n <= 100 then op(iquo(n, 2), [1/6, -1/30, 1/42, -1/30,
        5/66, -691/2730, 7/6, -3617/510, 43867/798, -174611/330,
        854513/138, -236364091/2730, 8553103/6, -23749461029/870,
        8615841276005/14322, -7709321041217/510, 2577687858367/6,
        -26315271553053477373/1919190, 2929993913841559/6,
        -261082718496449122051/13530, 1520097643918070802691/1806,
        -27833269579301024235023/690, 596451111593912163277961/282,
        -5609403368997817686249127547/46410,
        495057205241079648212477525/66,
        -801165718135489957347924991853/1590,
        29149963634884862421418123812691/798,
        -2479392929313226753685415739663229/870,
        84483613348880041862046775994036021/354,
        -1215233140483755572040304994079820246041491/56786730,
        12300585434086858541953039857403386151/6,
        -106783830147866529886385444979142647942017/510,
        1472600022126335654051619428551932342241899101/64722,
        -78773130858718728141909149208474606244347001/30,
        1505381347333367003803076567377857208511438160235/4686,
        -5827954961669944110438277244641067365282488301844260429/
        140100870,
        34152417289221168014330073731472635186688307783087/6,
        -24655088825935372707687196040585199904365267828865801/30,
        414846365575400828295179035549542073492199375372400483487/
        3318, -46037842994794576469355749690190468497942578727512\
        88919656867/230010, 1677014149185145836823154509786269900\
        207736027570253414881613/498, -20245761959352903602311311\
        60111731009989917391198090877281083932477/3404310, 660714\
        61941767865357384784742626149627783068665338893176199698\
        3/6, -131142648867401750799551142401931184334575027557202\
        8644296919890574047/61410, 117905727902108279988412335124\
        9215083775254949669647116231545215727922535/272118, -1295\
        58594820753752798942782853857674965934148371943514302331\
        6326829946247/1410, 1220813806579744469607301679413201203\
        958508415202696621436215105284649447/6, -2116004495972665\
        13097597728109824233673043954389060234150638733420050668\
        349987259/4501770, 67908260672905495624051117546403605607\
        342195728504487509073961249992947058239/6, -9459803781912\
        21252952274330694937218727028415330669361333856962043113\
        95415197247711/33330])
    else
        d := 4
             + trunc(evalhf((lnGAMMA(n + 1) - n*ln(2*Pi))/ln(10)))
             + length(n);
        lprint(`using ` . d . ` Digits`);
        s := trunc(evalhf(exp(0.5*d*ln(10)/n))) + 1;
        Digits := d;
        p := 1;
        t1 := 1.;
        t2 := t1;
        lprint(`start small prime loop at time` = time());
        while p <= s do
            p := nextprime(p);
            pn := p^n;
            pn1 := pn - 1;
            t1 := pn*t1;
            t2 := pn1*t2
        end do;
        gc();
        lprint(status);
        lprint(`used primes up to and including ` . p);
        lprint(`finish small prime loop at time` = time());
        z := t1/t2;
        gc();
        lprint(status);
        lprint(`finish full prec. division at time` = time());
        oz := 0;
        while oz <> z do
            oz := z;
            p := nextprime(p);
            Digits := max(d - ilog10(pn), 9);
            pn := Float(p,0);
            pn := p^n;
            pn1 := z/pn;
            Digits := d;
            z := z + pn1
        end do;
        gc();
        lprint(status);
        lprint(`used primes up to and including ` . p);
        lprint(`finish big prime loop at time` = time());
        p := evalf(2*pi);
        gc();
        lprint(status);
        lprint(`finish 2*Pi at time` = time());
        f := n!;
        gc();
        lprint(status);
        lprint(`finish factorial at time` = time());
        pn := p^n;
        gc();
        lprint(status);
        lprint(`finish (2*Pi)^n at time` = time());
        z := 2*z*f/pn;
        gc();
        lprint(status);
        lprint(
            `finish 2*z*n!/(2*Pi)^n (multiply and divide) at time`
             = time());
        s := 0;
        for p in numtheory[divisors](n) do
            p1 := p + 1; if isprime(p1) then s := s + 1/p1 end if
        end do;
        gc();
        lprint(status);
        lprint(`finish divisors of n loop at time` = time());
        s := frac(s);
        if irem(n, 4) = 0 then
            if s < 1/2 then z := -round(z) - s
            else z := -trunc(z) - s
            end if
        else
            s := 1 - s;
            if s < 1/2 then z := round(z) + s
            else z := trunc(z) + s
            end if
        end if;
        gc();
        lprint(status);
        lprint(`done at time` = time());
        z
    end if
end:

\end{verbatim}

\clearpage

\begin{table} \begin{center}
\begin{tabular}{|l|r|r|}
\hline
          Who               &    when    & highest $B_n$ \\
\hline
Bernoulli                   & 1713       &            10 \\
Euler                       & 1748       &            30 \\
J.C. Adams                  & 1878       &            62 \\
D.E. Knuth and Buckholtz    & 1967       &           360 \\
Greg Fee and Simon Plouffe  & 1996       &         10000 \\
Greg Fee and Simon Plouffe  & 1996       &         20000 \\
Greg Fee and Simon Plouffe  & 1996       &         30000 \\
Greg Fee and Simon Plouffe  & 1996       &         50000 \\
Greg Fee and Simon Plouffe  & 1996       &        100000 \\
Greg Fee and Simon Plouffe  & 1996       &        200000 \\
Simon Plouffe               & 2001       &        250000 \\
Simon Plouffe               & 2002       &        400000 \\
Simon Plouffe               & 2002       &        500000 \\
Simon Plouffe               & 2002       &        750000 \\
Berndt C. Kellner           & 2002       &       1000000 \\
Berndt C. Kellner           & 2003       &       2000000 \\
Pavlyk O.                   & 2005       &       5000000 \\
\hline
\end{tabular}
\caption{History of the computation of Bernoulli numbers}
\end{center} \end{table}

\begin{\bib}{99}

\bibitem{slo:ole} 

\begin{flushleft}
N. J. A. Sloane, ``The On-Line Encyclopedia of
Integer Sequences,'' available at
\href{http://www.research.att.com/~njas/sequence/}{\texttt http://www.research.att.com/$\sim$njas/sequences/}.
\end{flushleft}

\bibitem{A&S}
M.~Abramowitz and I.~Stegun, editors.\\
{\em Handbook of Mathematical Functions}. Dover, New York, 1970.

\bibitem{dilcher}
K. Dilcher, {\em, A Bibliography of Bernoulli Numbers},\\
\href{http://www.mscs.dal.ca/~dilcher/bernoulli.html}{\texttt http://www.mscs.dal.ca/~dilcher/bernoulli.html}.

\bibitem{GS} 
Xavier Gourdon et Patrick Sebah, {\em, Introduction  to Bernoulli Numbers},\\ \href{http://numbers.computation.free.fr/Constants/Miscellaneous/bernoulli.html}{\texttt http://numbers.computation.free.fr/Constants/Miscellaneous/bernoulli.html}.

\bibitem{hardy}
G. H. Hardy and E. L. Wright, {\em An Introduction to the Theory of
Numbers}, Oxford University Press, 1979.

\bibitem{kellner1)}
Kellner, B. C. {\em Über irreguläre Paare höherer Ordnungen}. Diplomarbeit. Göttingen, Germany: Mathematischen Institut der Georg August Universität zu Göttingen, 2002. \href{http://www.bernoulli.org/~bk/irrpairord.pdf}{\texttt http://www.bernoulli.org/~bk/irrpairord.pdf}.

\bibitem{kellner2)}
Kellner, B. C. {\em The Equivalence of Giuga's and Agoh's Conjectures}. 15 Sep 2004. \href{http://arxiv.org/abs/math.NT/0409259/}{\texttt http://arxiv.org/abs/math.NT/0409259/}. 

\bibitem{knuth} D. E. Knuth, {\em The Art of Computer Programming}
vol. 2, Addison-Wesley, Reading, MA, 1981 and 1997.

\bibitem{plouffe}
S. Plouffe, {\em The computation of big Bernoulli numbers}, from 1996 to 2002, on the author web site at \href{http://www.lacim.uqam.ca/$\sim$plouffe/Bigfiles/}{\texttt http://www.lacim.uqam.ca/$\sim$plouffe/Bigfiles/}.

\bibitem{rama}
S. Ramanujan, {\em Some properties of Bernoulli's numbers}, J. Indian Math. Soc., (1911), vol. 3, pp. 219-234.

\bibitem{Ribenboim}
P. Ribenboim, {\em The new Book of Prime Number Records}, Springer, (1996).

\bibitem{sloaneplouffe} 
Neil J.A. Sloane and Simon Plouffe,{\em, The Encyclopedia of Integer Sequences}, Academic Press, San Diego 1995, 587 pp. 

\bibitem{sloane} 
Neil Sloane and al. {\em The On-Line Encyclopedia of Integer Sequences}, available at \href{http://www.research.att.com/~njas/sequences/}{\texttt http://www.research.att.com/~njas/sequences/}.

\bibitem{weisstein1}
Eric Weisstein, {\em Agoh's conjecture},\\
\href{Mathworld, http://mathworld.wolfram.com/AgohsConjecture.html}{\texttt http://mathworld.wolfram.com/AgohsConjecture.html}.

\bibitem{weisstein2}
Eric Weisstein, {\em Bernoulli Numbers},\\
\href{Mathworld, http://mathworld.wolfram.com/BernoulliNumber.html}{\texttt http://mathworld.wolfram.com/BernoulliNumber.html}.

\end{\bib}

\bigskip
\hrule
\bigskip

\noindent 2000 {\it Mathematics Subject Classification}:
Primary 11B68; Secondary 05A10, 11A07, 11B64.

\noindent \emph{Keywords: } Bernoulli numbers, Euler product, Zeta function.

\bigskip
\hrule
\bigskip

\noindent (Concerned with sequences
\seqnum{A000367}, \seqnum{A000928}, \seqnum{A000928}, \seqnum{A002445}, \seqnum{A027641}, \seqnum{A027642})

\bigskip
\hrule
\bigskip

\vspace*{+.1in}
\noindent

\bigskip
\hrule
\bigskip

\noindent
Return to
\htmladdnormallink{Journal of Integer Sequences home page}{http://www.math.uwaterloo.ca/JIS/}.
\vskip .1in

\end{document}